\documentclass[11pt]{amsart}
\usepackage{amsmath,amsfonts,amsthm,amssymb,url,graphicx,algorithm}
\usepackage[noend]{algpseudocode}

\DeclareFontFamily{OT1}{rsfs}{}
\DeclareFontShape{OT1}{rsfs}{n}{it}{<-> rsfs10}{}
\DeclareMathAlphabet{\mathscr}{OT1}{rsfs}{n}{it}

\DeclareMathOperator{\num}{num}

\DeclareMathOperator{\den}{den}

\DeclareMathOperator{\mo}{\,mod}

\newtheorem*{main}{Main Theorem}

\newtheorem*{defn*}{Definition}

\numberwithin{equation}{section}
\title{An improved sieve of Eratosthenes}
\author{Harald Andr\'es Helfgott}
\address{Harald A. Helfgott, 
  Mathematisches Institut,
Georg-August Universit\"{a}t G\"{o}ttingen, Bunsenstra{\ss}e 3-5, D-37073 G\"{o}ttingen,
Germany; IMJ-PRG, UMR 7586,
  58 avenue de France, B\^{a}timent S. Germain, case 7012,
  75013 Paris CEDEX 13, France}
\email{harald.helfgott@gmail.com}
\begin{document}
\begin{abstract}
  We show how to carry out a sieve of Eratosthenes up to $N$ in space
  $O\left(N^{1/3} (\log N)^{2/3}\right)$ and time $O(N \log N)$.
  In comparison, the usual versions of the sieve
  take space about $O(\sqrt{N})$ and
  time at least linear on $N$. We can also apply our sieve to any subinterval
  of $\lbrack 1,N\rbrack$ of length
  $\Omega\left(N^{1/3}\right)$ in time close to linear on the length
  of the interval. Before, such a thing was possible only for subintervals of
  $\lbrack 1,N\rbrack$ of length $\Omega(\sqrt{N})$.
  
  Just as in (Galway, 2000), the approach here is
  related to Diophantine approximation, and also has close ties
  to Vorono\"i's work on the Dirichlet divisor problem. The advantage of the method here
  resides
  in the fact that, because the method we will give is based on the sieve of
  Eratosthenes, we will also be able to use it to factor integers,
  and not just to produce lists of consecutive primes.
\end{abstract}

\maketitle

\section{Introduction}

The sieve of Eratosthenes is a procedure for constructing all primes up
to $N$. More generally, such a sieve can be used
for factoring all integers up to $N$, or to compute the values $f(n)$,
$n\leq N$, of arithmetical functions $f$ that depend on the factorization
of integers. For instance, one can take $f=\mu$, the M\"obius function, or
$f=\lambda$, the Liouville function.

The point of the sieve of Eratosthenes is that it can be carried out in time
close to linear on $N$, even though determining whether an individual
integer is prime, let alone factoring it, takes much more than constant
time with current methods. Though a very na\"{\i}ve implementation of
the sieve would take space proportional to $N$, it is not hard to see how
to implement the sieve in space $O(\sqrt{N})$, simply by applying the
sieve to intervals of length $O(\sqrt{N})$ one at a time; the time taken
is still close to  linear\footnote{
  In its standard form,
  the sieve of Eratosthenes
  (segmented or not) takes time $O(N \log \log N)$; when used for sieving
  out primes, its time consumption can be reduced to $O(N)$.

  We follow the convention that
arithmetic operations (e.g., adding or multiplying two numbers) take
$O(1)$ time. This convention holds in the range in which using
a sieve is realistic. (Outside that range, such operations do take time
$O((\log N)^c)$, $c$ close to $1$.) We shall measure space in bits, again
reflecting how matters work in practice.}
  on $N$. (One could take shorter
intervals, but the algorithm would then become much less efficient.)
This is called the ``segmented sieve''; see \cite{1969-singleton-357} for an early reference.

Of course, the output still takes space linear on $N$,
but that is of less importance: we can store the output in slower memory
(such as a hard drive), or give the output in batches to a program that
needs it and can handle it sequentially.
Then the space used is certainly $O(\sqrt{N})$.



There has been a long series of improvements to the basic segmented sieve.
Most of them improve the running time or space by a constant factor
or by a factor in the order of $\log \log N$. Many work only when the
sieve is used to construct primes, as opposed to computing $\mu$, say.
See \cite{MR1726070}, which reviews the state of matters at the time
before improving on it; see also \cite{primesieve} for further references
and for a contemporary implementation combining some of the most useful
existing techniques.

In practice, saving space sometimes amounts to saving time, even when it seems, at
first, that there should be a trade-off.
As of the time of writing,
a good office computer can store about $10^{12}$ integers in very slow memory
(a hard drive), about  $10^9$ integers in memory working at intermediate speed
(RAM), and about $10^6$ integers in fast memory (cache); a program becomes
faster if it can run on cache, accessing RAM infrequently.
Having enough RAM is also an issue; sieves have been used, for instance, to
verify the binary Goldbach conjecture up to $4\cdot 10^{18}$ \cite{OSHP},
and so we are close to the point at which $O(\sqrt{N})$ space
might not fit in RAM.
Space constraints
can become more severe if several processor cores work in parallel and share resources.

Moreover, finding all primes within a short interval can be useful in itself.
For instance, there are applications in which we verify a conjecture on
one interval at a time, and we neither need nor can store a very long interval
in memory. See the verification of Goldbach's (binary) conjecture up to $4\cdot 10^{18}$ in \cite{OSHP}, and, in particular \cite[\S 1.2]{OSHP}.

Galway \cite{MR1850613} found a way to sieve using space $O(N^{1/3})$ and
time $O(N)$. Like the sieve in \cite{MR2031423}, on which
it is based\footnote{Atkin and Bernstein's preprint was already available in 1999, as the bibliography in \cite{MR1850613} states.}, Galway's sieve is specific to finding prime numbers. There is also the algorithm in
\cite{10.1007/11792086_15}, specific, again, to finding primes:
under the assumption of the Generalized Riemann Hypothesis, it
finds all primes up to $N$ in space $O((\log N)^3/\log \log N)$ and
time $O(N (\log N)^2/\log \log N)$; unconditionally, it runs in space
$O(N^{0.132})$ and time $O(N^{1.132})$. (It runs in time $O(N \log N)$ under the
assumption of a more specialized conjecture.)

We will show how to implement a sieve of Eratosthenes in space 
close to $N^{1/3}$ and still close to linear time. Our method is not limited
to finding prime numbers; it can be used to factor integers, or, of course,
to compute $\mu(n)$, $\lambda(n)$ or other functions given by the
factorization of $n$.

\begin{main}
  We can construct all primes $p\leq N$ in 
\begin{equation}\label{eq:rada1}
\text{space}\;\;O\left(N^{1/3} (\log N)^{2/3}\right)\;\;\;\;\;\;\;\;
\text{and}\;\;\;\;\;\; \text{time}\;\; O(N \log N).\end{equation}
We can also factor all integers $n\leq N$ in 
\begin{equation}\label{eq:rada2}
\text{space}\;\;O\left(N^{1/3} (\log N)^{5/3}\right)\;\;\;\;\;
\text{and}\;\;\;\;\;\; \text{time}\;\; O(N \log N).\end{equation}

  Moreover, for 
$N^{1/3} (\log N)^{2/3}
\leq \Delta\leq N$, we can construct all primes in an interval
  $\lbrack N-\Delta,N+\Delta\rbrack$ in
\[\text{space}\;\; O(\Delta)\;\;\;\;\;\text{and}\;\;\;\;\;\;
\text{time}\;\; O(\Delta \log N)\]
and factor all integers in the same interval in
\[\text{space}\;\; O(\Delta \log N)\;\;\;\;\;\;\text{and}\;\;\;\;\;\;
\text{time}\;\; O(\Delta \log N).\]
\end{main}
Here we recall that {\em space} refers to the number of bits used, and
{\em time} to the number of operations of words used (on integers of size
$O(N)$). 

The main ideas come from elementary number theory.
In order for us to be
able to apply the sieve to an interval $I$ of length $O(N^{1/3})$ without large
time inefficiencies, we need to be able to tell in advance which primes
(or integers) $d$
up to $\sqrt{N}$ divide at least one integer in $I$, without testing each $d$
individually. We can do this by Diophantine approximation, followed by
a local linear
approximation to the function $x\mapsto n/x$ for $n$ fixed, and
then by solving what amounts to a linear equation $\mo 1$.

The idea of using Diophantine approximation combined with a local linear approximation
is already present in \cite{MR2869058}, where
it  was used to compute $\sum_{n\leq x} \tau(n)$ in time
$O(x^{1/3} (\log x)^{O(1)})$. (We write
$\tau(n)$ for the number of divisors of an integer $n$.)
The basic underlying idea in
\cite{MR1850613} may be said to be the same as the one here:
we are speaking of a Diophantine
idea that stems ultimately from Vorono\"i's
work on the Dirichlet divisor problem \cite{zbMATH02656343} (in our work) and
Sierpinski's adaptation of the same method to the circle problem
\cite{sierpinski1906pewnem} (in the case of \cite{MR1850613}).\footnote{
To be precise, the immediate inspiration for \cite{MR2869058} came
from Vinogradov's simplified version of Vorono\"i's method, as in
\cite[Ch. III, exercises 3--6]{MR0062138}. A 
bound of roughly the same quality as Vorono\"i's result was claimed long
before Vorono\"i in
\cite{pfeiffer}.
 While the proof there was apparently incomplete, Landau
later showed \cite{zbMATH02625237} that it could be made into an actual proof,
and sharpened to give a result matching Vorono\"i's.
Thanks are due to S.\ Patterson for pointing out Pfeiffer's result to the 
author.}
For that matter,  \cite[\S 5]{MR1850613} already suggests that Vorono\"i's
work on the Dirichlet divisor problem could be used to make the
sieve of Eratosthenes in space about $O(N^{1/3})$ and close to linear time.
We should also make clear that Galway can sieve out efficiently
segments of length about $x^{1/3}$, just as we do.

One difference between this paper and Galway's
is that the relation to Vorono\"i's and Sierpinski's work in
Galway's paper may be said to be
more direct, in that Galway literally dissects a region
  between two circles, much as Sierpinski does. In the case of the
  present paper, we can say that Vorono\"i's 
 main idea originated in the context of giving
 elementary estimates (for $\sum_{n\leq x} \tau(n)$)
 and is now used to carry out
  an exact computation.

  Another precedent that must be mentioned is Oliveira e Silva's technique for
  efficient cache usage \cite{OeS}, \cite[Algorithm 1.2]{OSHP}.
    It seems that
  this technique can be combined with the algorithm here. Such a
  combination can be useful, for instance, when $N$ is so large that $N^{1/3}$ bits fit in RAM but not in cache.

  {\bf Additional motivation.} My initial interest
  in the problem stemmed from the fact that I had to compute values of $\mu(n)$
  so as to check inequalities of the form\footnote{Mertens's conjecture
    states that the inequality $\sum_{n\leq x} \mu(n) \leq \sqrt{x}$ 
    holds for all $x$. That conjecture has been
    disproved \cite{MR783538}, but the inequality is known to hold for all
    $x\leq 10^{16}$
    \cite{Hurst}, and may hold in a far wider range.}
    $\sum_{n\leq x} \mu(n) \leq \sqrt{x}$,
    $\sum_{n\leq x} \mu(n)/n \leq 2/\sqrt{x}$, etc., for all $x$ less
  than some large finite $N$. The importance of such sums in number theory is
  clear. Explicit results on them for $x$ bounded help complement analytic
  estimates, which are generally strong only when $x$ is large.

  While we can obviously determine values of $\mu(n)$ by first factorizing
  $n$ and then using the definition of $\mu$, it is also possible and
  rather simple
  to save some space and time in practice by modifying the procedure
  we will give (Algorithm \ref{alg:segsievefac})
  so as to keep track of $\mu(n)$ instead of the list of factors.
  Time and space complexity remain the same.
  
  {\bf Conventions and notation.}
  Integer operations are assumed to take constant time.
  As is usual,
  we write $f(x) = O(g(x))$ to mean that there exists a constant $C>0$
  such that $|f(x)|\leq C g(x)$ for all large enough $x$. We write
  $f(x) = O^*(g(x))$ to mean that $|f(x)|\leq g(x)$ for all $x$.
  We use either $f(x) \ll g(x)$ or $g(x) \gg f(x)$ to mean $f(x) = O(g(x))$.
    
  For $n$ a non-zero integer and $p$ a prime, $v_p(n)$ denotes the largest
  $k$ such that $p^k|n$. As is customary,
  we write $(a,b)$ for the gcd (greatest common divisor)
  of $a$ and $b$, provided that no confusion with the ordered pair $(a,b)$
  is possible.

  Given $\alpha\in \mathbb{R}$, we denote by
 $\{\alpha\}$ the element of $\lbrack 0,1)$ congruent to
$\alpha$ modulo $1$, i.e., modulo $\mathbb{Z}$.

  {\bf Acknowledgements.}
  The author is currently supported by funds from his Humboldt Professorship.
Thanks are due to Manuel Drehwald, for having written code implementing
an earlier version
of the algorithm, to Lola Thompson, for helpful suggestions, and to
an anonymous referee, for several useful remarks.

\section{Analysis of the problem}\label{sec:analys}

Let $I = \lbrack n-\Delta, n+\Delta\rbrack \subset \lbrack 0,x\rbrack$ be
an interval. How can we tell which integers
$m\leq \sqrt{x}$ have at least one integer multiple in the interval $I$?

Our motivation for asking this question is that we will be taking intervals $I$
and sieving them by all integers, or all primes, $m\leq \sqrt{x}$.
Going over all integers $\leq \sqrt{x}$ would take time at least $\sqrt{x}$,
which could be much larger than $\Delta$.


We will be able to sieve our intervals by $m\leq \sqrt{x}$ by producing
a list of those $m$ which might have multiples in $I$,
without testing all $m\leq \sqrt{x}$ in succession. A quick probabilistic heuristic
hints that the number of such $m$ should be proportional to $\Delta \log x$,
which, for $\Delta\sim x^{1/3}$, is much smaller than $\sqrt{x}$.

Let $K\geq 2$ be a parameter to be set later.
If $m\leq K \Delta$,  
we simply sieve by $m$. Doing so takes time $O(1+\Delta/m)$. We could,
of course, test first whether $m$ is prime; we will discuss this option later.

Assume henceforth that $m>K \Delta$.
The interval $I$ contains a multiple of $m$ if and only if
\[\left\{\frac{n}{m}\right\} \in \left\lbrack -\frac{\Delta}{m},\frac{\Delta}{m}\right\rbrack \mo 1,\]
that is, $\{n/m\} \in \lbrack -\Delta/m,\Delta/m\rbrack + \mathbb{Z}$.

  Say we have already dealt with all $m\leq M$, where $M\geq K \Delta$,
  and that we want to examine $m$ close to $m_0$, where $m_0>M$.
Write $m = m_0 + r$.
The truncated Taylor expansion
    \[f(m) = f(m_0) + f'(m_0) r + \frac{f''(m_0 + \theta r)}{2} r^2\]
    (where $\theta$ is some element of $\lbrack 0,1\rbrack$) gives us, once
    we set $f(x) = n/x$,
    \begin{equation}\label{eq:anzio}
      \frac{n}{m} = \frac{n}{m_0} - \frac{n}{m_0^2} r + O^*\left(\frac{n}{
        m_-^3} r^2\right),
    \end{equation}
    where $m_- = \min(m,m_0)$.
    We will make sure that $r$ is small enough that
    $n r^2/m_-^3 \lesssim \kappa \Delta/m$,
    where $\kappa>0$.
    (That is, we allow our error term to be not much larger than the interval
    we are trying to hit.) 
    We ensure $r$ satisfies this condition by letting
    \begin{equation}\label{eq:shapewater}
      R = \left\lfloor \sqrt{\frac{\kappa \Delta}{n}} M\right\rfloor,\;\;\;\;\;\;
    m_0 = M + R.\end{equation}
    Then $n r^2/m_-^3 \leq n r^2/M^3 \leq \kappa \Delta/M$ for all $m=m_0+r$,
    $r\in \lbrack -R,R\rbrack$.
    We let $\kappa=1/4$, since it is a value that is neither too large nor too
    small. We will also assume $\Delta \geq n^{1/3}$, and so
    $R\geq \lfloor K/2\rfloor \geq 1$, since $K\geq 2$.

    
    We have thus reduced our problem to that
    of quickly finding all $r$ in an interval $\lbrack -R,R\rbrack$ such
    that $P(r) \in \lbrack-\eta,\eta\rbrack \mo 1$, where
    $P(r)$ is the linear polynomial $- (n/m_0^2) r + (n/m_0)$ and
    $\eta = (1+\kappa) \Delta/M = 5 \Delta/4 M$.
    In other words, we are being asked to find all approximate solutions
    to a linear equation in $\mathbb{R}/\mathbb{Z}$ efficiently. 
    Let us
    assume that $K\geq 5/2$, so that $\eta\leq 1/2$; otherwise there is not much
    to do.
    
    Let $\alpha_1 = - (n/m_0^2)$, $\alpha_0 = n/m_0$.
  Given a rational number
  $\alpha$, we can find -- by means of continued fractions --
  an approximation $a/q$ to $\alpha$ with $q\leq Q$ and
  \begin{equation}\label{eq:jukuru}
    \left|\frac{a}{q} - \alpha\right|\leq \frac{1}{q Q}\end{equation}
  in time $O(\log Q)$.
  The procedure {\bf DiophAppr}$(\alpha,Q)$ to obtain $a/q$ is given in Algorithm \ref{alg:diophfrac}; it runs in time $O(\log Q)$ and constant space.
  The fact that its output satisfies (\ref{eq:jukuru})  follows from \cite[Thm. 164]{zbMATH03657869}
  or \cite[Thm. 9]{MR1451873}.
  (In the notation of both sources:
      since any two consecutive approximants $p_n/q_n$, $p_{n+1}/q_{n+1}$
      to $\alpha$ satisfy $|p_n/q_n - \alpha|\leq 1/q_n q_{n+1}$,
  the last $p_n/q_n$ with $q_n\leq Q$ satisfies $|p_n/q_n-\alpha|<1/q_n Q$.)

  We invoke {\bf DiophAppr}$(\alpha_1,2R)$, and obtain a rational
  $a/q$ with $(a,q)=1$ and $q\leq Q = 2 R$ satisfying (\ref{eq:jukuru}) for $\alpha = \alpha_1$.
  Hence, for $r\in \lbrack - R,R\rbrack$,
  \begin{equation}\label{eq:garla}
    P(r) = \alpha_1 r + \alpha_0 = \alpha_0 + \frac{a r}{q} + O^*\left(\frac{1}{2q}\right) \equiv \frac{c + a r}{q} + O^*\left(\frac{1}{q}\right) \mo 1
    \end{equation}
  for $c = \lfloor \alpha_0 q + 1/2\rfloor$.
  We have thus reduced our problem to a problem in $\mathbb{Z}/q\mathbb{Z}$:
  if $P(r) \in \lbrack-\eta,\eta\rbrack$, then
  \begin{equation}\label{eq:norar}
    c + a r \in \{-k-1,-k,\dotsc,k+1\} \mo q,\end{equation}
  where $k = \lfloor\eta q\rfloor$. We must find the values of $r$
  satisfying (\ref{eq:norar}).

  The solutions to (\ref{eq:norar}) are clearly given by
  \begin{equation}\label{eq:amark}
    r \equiv - a^{-1} (c + j) \mo q
  \end{equation}
  for $-k-1\leq j\leq k+1$. The multiplicative inverse $a^{-1} \mo q$
  of $a$ can be computed in time $O(\log q)$ by the Euclidean algorithm.
In fact, we compute it
 at the same time as the continued fraction that gives us $a/q$: by \cite[Thm. 2]{MR1451873}, two consecutive approximants
    $p_n/q_n$, $p_{n+1}/q_{n+1}$ satisfy $p_n^{-1} = (-1)^{n+1} q_{n+1} \mo q_{n}$;
moreover, $p_n^{-1} = (-1)^{n+1} q_{n-1} \mo q_n$. 

Thus, we simply need to go over all $m$ of the form $m_0+r$, where
$r$ goes over integers in $\lbrack -R,R\rbrack$ satisfying (\ref{eq:amark}).
Since $m_0-R = M > K \Delta \geq 2 \Delta$, each such $m$ has at most
one multiple in $I=\lbrack n-\Delta,n+\Delta\rbrack$. We do not miss any
$m\in \lbrack m_0-R,m_0+R\rbrack = \lbrack M,M+2 R\rbrack$ that has a
multiple in $I$.

Due to the errors involved in the truncation of the Taylor expansion 
and in Diophantine approximation, we may obtain some $m$ that have no
multiples within $I$. We simply ignore them, after checking that that is the
case.

\section{Description of the algorithm}




We have already given a nearly full description of the algorithm. It only
remains to give the pseudocode (Algorithms \ref{alg:segsieve}--\ref{alg:diophfrac}), with some commentary.


\begin{algorithm}
  \caption{Main algorithm: 
sieving $\lbrack n-\Delta,n+\Delta\rbrack \subset \mathbb{R}^+$}\label{alg:segsieve}
  \begin{algorithmic}[1]
\Function{NewSegSiev}{$n$,$\Delta$,$K$}
    \Ensure{$S_j=1$ if $n+j$ is prime, $S_j=0$ otherwise, for $-\Delta\leq j
\leq \Delta$}
    \Require{$n,\Delta\in \mathbb{Z}^+$,   $\sqrt[3]{n}\leq \Delta<n$,
      $K\geq 5/2$}
\State{$S' \gets \text{\textsc{SubSegSiev}}(n-\Delta,2 \Delta,K \Delta)$}\Comment{sieve by all $p\leq K\Delta$}
\State{$S_j \gets S'_{j+\Delta}$ for all $-\Delta \leq j\leq \Delta$}
\State{$M\gets \lfloor K \Delta \rfloor + 1$}
\While{$M\leq \sqrt{n+\Delta}$}
\State{$R\gets \lfloor M \sqrt{\Delta/ 4 n}\rfloor$, $m_0 \gets M+R$}
\State{$\alpha_1 \gets \{- n/m_0^2\}$, $\alpha_0 \gets \{n/m_0\}$, 
  $\eta\gets 5 \Delta/4 M$}
\State{$(a,a^{-1},q) \gets \text{\textsc{DiophAppr}}(\alpha_1,2 R)$}
\State{$c\gets \lfloor \alpha_0 q + 1/2\rfloor$,
  $k\gets \lfloor \eta q\rfloor$}
\For{$-k-1\leq j\leq k+1$}
\State{$r_0\gets - a^{-1} (c+j) \mo q$}
\For{$m\in (m_0 + r_0 + q \mathbb{Z}) \cap \lbrack M,M+2R\rbrack$}
\State{$n'\gets \lfloor (n+\Delta)/m\rfloor\cdot m$} \Comment{$n'$ is a multiple of $m$}
\If{$n' \in \lbrack n-\Delta,n+\Delta\rbrack \wedge (n'>m)$}
\State{$S_{n'-n_0}\gets 0$} \Comment{thus sieving out $n'$}
\EndIf
\EndFor
\EndFor
\State{$M\gets M+2 R+1$}
\EndWhile
\State{\Return{$S$}}
\EndFunction
\end{algorithmic}
\end{algorithm}

\begin{algorithm}
  \caption{A simple sieve of Eratosthenes}\label{alg:erasieve}
  \begin{algorithmic}[1]
    \Function{SimpleSiev}{$N$} 
    \Ensure{for $1\leq n\leq N$, $P_n=1$ if $n$ is prime, $P_n=0$ otherwise}
    \State{$P_1\gets 0$, $P_2\gets 1$, $P_n\gets 0$ for $n\geq 2$ even,
      $P_n\gets 1$ for $n\geq 3$ odd}
    \State{$m\gets 3$, $n\gets m\cdot m$}
    \While{$n\leq N$}
    \If{$P_m=1$}
    \While{$n\leq N$} \Comment{[sic]}
    \State{$P_n\gets 0$, $n\gets n + 2 m$} \Comment{sieves odd multiples
    $\geq m^2$ of $m$}
    \EndWhile
    \EndIf
    \State{$m\gets m+2$, $n\gets m\cdot m$}
    \EndWhile
    \State{\Return{$P$}}
    \EndFunction
    \vskip 3pt
    \noindent {\bf Time:} $O(N \log \log N)$.\; {\bf Space:} $O(N)$.
  \end{algorithmic}
\end{algorithm}

\begin{algorithm}
   \caption{Segmented sieves of Eratosthenes, traditional version}\label{alg:oldsegsieve}
  \begin{algorithmic}[1]
\Function{SimpleSegSiev}{$n$,$\Delta$,$M$}
\Comment{sieves $\lbrack n,n+\Delta\rbrack$ by primes $p\leq M$}
\Ensure{ $S_j=\begin{cases} 0 &\text{if $n+j=0,1$ or if $p|(n+j)$ for some $p\leq M$}\\ 1 &\text{otherwise}\end{cases}$}
\State{$S_j \gets 1$ for all $0\leq j\leq \Delta$}
\State{$S_j\gets 0$ for $0\leq j\leq 1-n$} \Comment{[sic; excluding $0$ and $1$ from prime list]}
       \State{$P\gets \text{\textsc{SimpleSiev}}(M)$}
    \For{$1\leq m\leq M$}
    \If{$P_m=1$}
    \State{$n' \gets \max(m\cdot \lceil n/m\rceil,2 m)$}
    \While{$n' \leq n+\Delta$}\Comment{$n'$ goes over mults.\ of $m$ in $n+\lbrack 0,\Delta\rbrack$}
\State{$S_{n'-n} \gets 0$, $n' \gets n' + m$}
\EndWhile
\EndIf
\EndFor
\State{\Return{$S$}}
\EndFunction
\vskip 3pt
\noindent {\bf Time:} $O((M + \Delta) \log \log M)$.\; {\bf Space:}
$O(M + \Delta)$.

\vskip 5pt

\Function{SubSegSiev}{$n$,$\Delta$,$M$}
\Comment{sieves $\lbrack n,n+\Delta\rbrack$ by primes $p\leq M$}
\Ensure{ $S_j=\begin{cases} 0 &\text{if $n+j=0,1$ or if $p|(n+j)$ for some $p\leq M$}\\ 1 &\text{otherwise}\end{cases}$}
\State{$S_j \gets 1$ for all $\max(0,2-n)\leq j\leq \Delta$}
\State{$S_j\gets 0$ for $0\leq j\leq 1-n$}
\State{$\Delta' \gets \lfloor \sqrt{M}\rfloor$, $M'\gets 1$}
\While{$M'\leq M$}
\State{$P\gets \text{\textsc{SimpleSegSiev}}(M',\Delta',\lfloor \sqrt{M'+\Delta'}\rfloor)$}
\For{$M'\leq p< \min(M,M'+\Delta')$}
\If{$P_{p-M'}=1$} \Comment{if $m$ is a prime\dots}
\State{$n' \gets \max(p\cdot \lceil n/p\rceil, 2 p)$} 
\While{$n' \leq n+\Delta$} 
\State{$S_{n'-n} \gets 0$, $n' \gets n' + p$}
\EndWhile
\EndIf
\EndFor
\State{$M'\gets M'+\Delta'+1$}
\EndWhile
\State{\Return{$S$}}
\EndFunction
\vskip 3pt
\noindent {\bf Time:} $O((M+\Delta) \log \log M)$.\; {\bf Space:} $O(\sqrt{M} + \Delta)$.

\vskip 5pt

\Function{SegSiev}{$n$,$\Delta$} \Comment{finds primes in $\lbrack n, n+\Delta\rbrack$}
\Ensure{for $0\leq j\leq \Delta$,
  $S_j=1$ if $n+j$ prime, $S_j=0$ otherwise}
\State{\Return{\textsc{SubSegSiev}$(n,\Delta,\lfloor \sqrt{n+\Delta}\rfloor)$}}
\EndFunction
\vskip 3pt
\noindent {\bf Time:} $O((\sqrt{n}+\Delta) \log \log (n+\Delta))$.
\Comment $\sqrt{n+\Delta} \leq \sqrt{n}+\sqrt{\Delta}\leq \sqrt{n}+\Delta$
\vskip 1pt
\noindent {\bf Space:} $O(n^{1/4} + \Delta)$.
\Comment $(n+\Delta)^{1/4} \leq n^{1/4}+\Delta$
\end{algorithmic}
\end{algorithm}

\begin{algorithm}
  \caption{Finding a Diophantine approximation via continued fractions}\label{alg:diophfrac}
  \begin{algorithmic}[1]
    \Function{DiophAppr}{$\alpha$,$Q$}
    \Ensure{Returns $(a,a^{-1},q)$ s.t.
      $\left|\alpha-\frac{a}{q}\right|\leq \frac{1}{q Q}$, $(a,q)=1$, $q\leq Q$, $a a^{-1} \equiv 1 \mo q$}
    \State{$b\gets \lfloor \alpha\rfloor$, $p\gets b$, $q\gets 1$,
    $p_-\gets 1$, $q_-\gets 0$, $s\gets 1$}
    \While{$q\leq Q$}
    \If{$\alpha=b$}
    \Return{$(p,- s q_-, q)$}
    \EndIf
    \State{$\alpha \gets 1/(\alpha-b)$}
    \State{$b\gets \lfloor \alpha\rfloor$,
      $(p_+,q_+) \gets b\cdot (p,q) + (p_-,q_-)$}
    \State{$(p_-,q_-)\gets (p,q)$, $(p,q)\gets (p_+,q_+)$, $s\gets -s$}
    \EndWhile
    \State{\Return{$(p_-, s q,q_-)$}}
    \EndFunction
    \vskip 3pt
\noindent {\bf Time:} $O(\log \max(Q,\den(\alpha)))$.\; {\bf Space:} $O(1)$.
      \end{algorithmic}
\end{algorithm}

All variables are integers, rationals, or tuples or sets with
integer or rational entries. Thus, all arithmetic operations can be
carried out exactly. We write $\num(\alpha)$ and $\den(\alpha)$ for
the numerator and denominator of $\alpha \in \mathbb{Q}$ (written
in minimal terms: $\num(\alpha)$, $\den(\alpha)$ coprime, $\den(\alpha)>0$).

In Algorithm \ref{alg:segsievefac}, 
we are defining $\alpha_1 \gets \{- n/m_0^2\}$, $\alpha_0 \gets \{n/m_0\}$,
rather than $\alpha_1 \gets - n/m_0^2$, $\alpha_0 \gets n/m_0$, as in the
exposition above.
The motivation is simply to keep variable sizes small. It is easy to see that
the values of $\alpha_0, \alpha_1$ matter only $\mo \mathbb{Z}$.

\subsection{Sieving for primes}

The main function is \textsc{NewSegSiev}
(Algorithm \ref{alg:segsieve}).
It takes as inputs $n$ and $\Delta$, as well as the parameter $K$, which
affects time and space consumption.
The basic procedure is easy to summarize.
\textsc{NewSegSiev} uses \textsc{SubSegSiev} (Algorithm \ref{alg:oldsegsieve}),
a segmented sieve of traditional type, to sieve for primes $p\leq K \Delta$. 
It then proceeds as detailed in \S \ref{sec:analys} to find the integers $m$
in $\lbrack M, M+2 R\rbrack$ that may divide an integer
in the interval $\lbrack n-\Delta,n+\Delta\rbrack$. Then it sieves by
them. The key step, in finding such $m$,
is to call \textsc{DiophAppr}, which
uses continued fractions in a standard way to find (i)
a rational approximation $a/q$ to the input $\alpha_1$, (ii)
the inverse $a^{-1} \mo q$.

We could try to improve on \textsc{NewSegSieve} 
by throwing out even values of $m$, say; as they are certainly not prime,
we need not sieve by them. More details are given in ``Side notes on wheels'',
below.

{\em Preexistent sieves as subroutines.} Going back to \textsc{SubSegSiev}: we could avoid using it at all,
just by sieving by all integers $m\leq K \Delta$, rather
than by all primes $p\leq K \Delta$. That would give a running time
of $O(K \Delta + \Delta \log M)$ rather than $O(K \Delta \log \log M)$.
We take the slightly more complicated route in Algorithm \ref{alg:oldsegsieve},
not just because it is better for $K = o((\log M)/\log \log M)$,
but also for the sake of exposition, in that we get to see several existing
forms of the sieve of Eratosthenes. Note, however, that none of this will
decrease the order of magnitude of the time taken by our entire algorithm,
since the total time will be at least in the order of $\Delta \log M$.

Function \textsc{SimpleSiev} is a relatively simple kind of sieve of
Eratosthenes. It sieves the integers up to $N$ by the primes up to
$\sqrt{N}$, where these primes are found by this very same process.
It is clear that we need to sieve only by the primes up to $\sqrt{N}$,
since any composite number $n\leq N$ has at least one prime factor
$p\leq \sqrt{N}$.
Sieving only by the primes, rather than by all integers,
is enough to take down the running time to $O(N \log \log N)$.
We also use a very primitive ``wheel'' in that we sieve using only
odd multiples of the primes. (Again, see the comments on wheels below.)

The basic segmented sieve is implemented as
\textsc{SimpleSegSiev}$(n,\Delta,M)$.
It uses \textsc{SimpleSiev} so as to determine the primes up to $M$.
Then it sieves the interval $\lbrack n,n+\Delta\rbrack$ by them.

We could use \textsc{SimpleSegSiev} instead of \textsc{SubSegSiev}.
The point of \textsc{SubSegSiev} is simply to reduce space consumption
by one more iteration: \textsc{SubSegSiev} determines the primes up to $M$ by
\textsc{SimpleSegSiev}, taking only space $O(\sqrt{M})$ at a time; it sieves
$\lbrack n,n+\Delta\rbrack$
by these primes as it goes along. The total time
taken by calls on \textsc{SimpleSegSiev} is $O(M \log \log M)$; to this
we add the time $O(\Delta \log \log M)$
taken by sieving the interval $\lbrack n,n+\Delta\rbrack$ by the primes up to $M$.

Incidentally: it should be clear that all factors of $\log \log N$, $\log \log M$ and the like are coming from Mertens's classical asymptotic statement
\[\sum_{p\leq N} \frac{1}{p} = \log \log N + O(1).\]
For instance, the number of times the
instructions $S_{n'-n} \gets 0$, $n' \gets n' + m$
are executed in \textsc{SimpleSegSiev} or \textsc{SubSegSiev}
is at most
\[\begin{aligned}\sum_{p\leq M} \left\lceil \frac{\Delta+1}{p}\right\rceil &\leq
\Delta \sum_{p\leq M} \frac{1}{p} + \sum_{p\leq M} 1 \\
&= O(\Delta \log \log M + M).\end{aligned}\] 

{\em Side note on wheels.} In general, a ``wheel'' is just
$(\mathbb{Z}/P \mathbb{Z})^*$, where $P = \prod_{p\leq C} p$ for some constant
$C$. We would use it by sieving only by multiples $m\cdot p'$ of the
primes $p'$,
where $m$ reduces $\mo P$ to an element of the wheel. Obviously $m \mo P$
should be constantly updated by shifting as $m$ increases, rather than be
determined by division each time; hence the ``wheel''. It is possible (and
common, at least in theoretical analyses) to choose $C$ of size
$\delta \log N$ for a small constant $\delta$
(say, $\delta=1/4$), with the consequence that $P\sim N^\delta$.
Then only a proportion $\prod_{p\leq C} (1-1/p) \sim \frac{e^{-\gamma}}{\log C}
= O(1/\log \log N)$ of the elements of
$\mathbb{Z}/P\mathbb{Z}$ (here $\gamma$ is Euler's constant)
are in $(\mathbb{Z}/P \mathbb{Z})^*$. With appropriate coding, this fact can
be used to reduce the running time of a simple sieve or a segmented
sieve for primes (such as \textsc{SimpleSiev} or \textsc{SimpleSegSiev})
by a factor of $\log \log N$ \cite{pritchard1983fast}.
We will not bother including wheels in our pseudocode, since they would
reduce only a second-order term of the total time bound in this fashion; the 
time complexity of \textsc{NewSegSiev} would remain the same. They
can obviously be added in implementation.

It is tempting to introduce a wheel in a different place, namely,
to make sure that the variable $m$
in \textsc{NewSegSiev} has no prime factors
$p<\delta \log N$. The hope there would be to reduce the main term in the
total time by a factor of $\log \log N$. However, as we will later see,
$q$ will be usually close to $2 R$; thus, most of the time, it will not
make sense to use a large wheel as $m$ goes over integers in
$(m_0+r_0 + q\mathbb{Z})\cap \lbrack M,M+2R\rbrack$, as there will be
few such integers. Using a wheel on $q$ instead -- something that would
take time $\gg q$ to set up -- would make no sense:
we do not go over $\mathbb{Z}/q\mathbb{Z}$ several times, but rather
go over a few elements of it once.

What can make sense is introduce a very small wheel, of bounded size, to
attempt to gain a constant factor. For instance, a wheel of size $2$
would work as follows: if $q$ is even, then, in the loop on $j$
in \textsc{NewSegSiev}, we consider only values of $j$ such that
$c+j$ is odd; if $q$ is odd, then we hop over every other value of $m$
in a given congruence class mod $q$ within $\lbrack M,M+2R\rbrack$,
considering only odd values $m$. This sort of modification will reduce total
time consumption only by a constant factor, and so, for the sake of
simplicity, we do not include it in the pseudocode, or use it further.

{\em Variable size.} We should bound the size of our variables in case
we choose to implement our algorithm using fixed-size integers and rationals. It is easy
to see that our integers will be of size $\leq n+\Delta$. We should
also bound the size of our rational variables.
It is trivial to modify function \textsc{NewSegSiev}
so that $m_0$ is always $\leq \sqrt{n+\Delta}$.
It is then easy to show that
we work entirely with integers -- in numerators, in denominators
or on their own -- of size
$\leq \min(n+\Delta,(2 R)^2)\leq n+\Delta$.

For $n$ very large and $\Delta$ not as large, it may be helpful to
store some variables (such as $a$, $a^{-1}$, $q$, $k$, $j$ and $r_0$)
in smaller integer types
($64$ bits, say); these variables are all bounded by
$2 R \leq M \sqrt{\Delta/n}\leq \sqrt{\Delta (n+\Delta)/n}$.



\subsection{Sieving for factorization}

We will now see how to modify our algorithm so that it factorizes all integers
in the interval $\lbrack n-\Delta,n+\Delta\rbrack$, rather than simply finding all
primes in that interval. Time and space usage will not be much greater
than when we are just finding primes.
It goes without saying that this makes it possible
to compute various arithmetic functions (the M\"obius function $\mu(m)$,
the Liouville function $\Lambda(m)$, etc.) for all $m\in \lbrack n-\Delta,n+\Delta\rbrack$.

For the sake of clarity, we first give a well-known procedure for
factoring all integers in an interval $\lbrack n,n+\Delta\rbrack$ by
means of the sieve of Eratosthenes (Algorithm \ref{alg:basterno}), just
as we went over a traditional segmented sieve (Algorithm \ref{alg:oldsegsieve})
before describing our sieve for primes.
We will later reuse most of the subroutines.

\begin{algorithm}
   \caption{Segmented sieve of Eratosthenes for factorization (traditional)}\label{alg:basterno}
  \begin{algorithmic}[1]
\Function{SubSegSievFac}{$n$,$\Delta$,$M$}
\Comment{finds prime factors $p\leq M$}
\Ensure{for $0\leq j\leq \Delta$, $F_j=\{(p,v_p(n+j))\}_{p\leq M, p|n+j}$}
\Ensure{for $0\leq j\leq \Delta$, $\Pi_j = \prod_{p\leq M, p|(n+j)} p^{v_p(n+j)}$.}
\State{$F_j \gets\emptyset$, $\Pi_j\gets 1$ for all $0\leq j\leq \Delta$}
\State{$\Delta'\gets \lfloor \sqrt{M}\rfloor$, $M'\gets 1$}
\While{$M'\leq M$}
\State{$P\gets \text{\textsc{SimpleSegSiev}}(M',\Delta',\lfloor \sqrt{M'+\Delta'}
  \rfloor)$}
\For{$M'\leq p<M'+\Delta'$}
\If{$P_{p-M'}=1$}  \Comment{if $p$ is a prime\dots}
\State{$k\gets 1$, $d\gets p$} \Comment{$d$ will go over the powers $p^k$ of $p$}
\While{$d\leq n+\Delta$}
\State{$n'\gets d\cdot \lceil n/d\rceil$}
\While{$n'\leq n+\Delta$}
\If{$k=1$}
\State{{\bf append} $(p,1)$ to $F_{n'-n}$}
\Else
\State{{\bf replace} $(p,k-1)$ by $(p,k)$ in $F_{n'-n}$}
\EndIf
\State{$\Pi_{n'-n}\gets p\cdot \Pi_{n'-n}$, $n'\gets n'+d$}
\EndWhile
\State{$k\gets k+1$, $d\gets p\cdot d$}
\EndWhile
\EndIf
\EndFor
\State{$M'\gets M'+\Delta'$}
\EndWhile
\State{\Return{$(F,\Pi)$}}
\EndFunction
\vskip 3pt
\noindent {\bf Time:} $O((M + \Delta) \log \log (n+\Delta))$,
\vskip 1pt
\noindent {\bf Space:}
$O(M + \Delta \log (n+\Delta))$.
\vskip 5pt

\Function{SegSievFac}{$n$,$\Delta$}
\Comment{factorizes all $n'\in \lbrack n, n+\Delta\rbrack$}
\Ensure{for $0\leq j\leq \Delta$, $F_j$ is the list of pairs $(p,v_p(n+j))$
  for $p|n+j$}
\State{$(F,\Pi) \gets \textsc{SubSegSievFac}(n,\Delta,\lfloor \sqrt{n+\Delta}
  \rfloor)$}
\For{$n\leq n'\leq n+\Delta$}
\If{$\Pi_{n'-n} \ne n'$}
\State{$p_0\gets n'/\Pi_{n'-n}$, {\bf append} $(p_0,1)$ to $F_{n'-n}$}
\EndIf
\EndFor
\State{\Return{$F$}}
\EndFunction
\vskip 3pt
\noindent {\bf Time:} $O((\sqrt{n} + \Delta) \log \log (n+\Delta))$, 
\vskip 1pt
\noindent {\bf Space:}
$O(n^{1/4} + \Delta \log (n+\Delta))$.
  \end{algorithmic}
\end{algorithm}

Our new sieve for factoring
(\textsc{NewSegSievFac},
Algorithm \ref{alg:segsievefac}), designed for intervals around $x$
of length $\Delta\gg x^{1/3}$, is very similar
to our sieve for primes (Algorithm \ref{alg:segsieve}).
We use a classical segmented
sieve (\textsc{SubSegSievFac}, Algorithm \ref{alg:basterno})
to find all factors $p\leq K\Delta$ of $n+j$ for $-\Delta\leq j\leq \Delta$.
After the call to $\textsc{SubSegSievFac}$, the variable $\Pi_j$ contains
$\prod_{p\leq K \Delta} p^{v_p(n+j)}$.
Since $K \Delta > (2 n)^{1/3} > (n+\Delta)^{1/3}$,
we see that $(n+j)/\Pi_j$ is either $1$ or the product of at most two
primes $> K \Delta$ (not necessarily distinct).
In the innermost loop of \textsc{SubSegSievFac},
when we come across an $m>K\Delta$ that divides not just $n+j$, but
$(n+j)/\Pi_j$,
we multiply $\Pi_j$ by $m$, or by its square, if $m^2|(n+j)/\Pi_j$,
and include $m$ (or its square) in the factorization.
Note that $m$
has to be a prime, or else $\Pi_j$ would have already been
multiplied by some factor $p^{v_p(n+j)}$, $p|m$, $p<m$, either earlier in the loop
or in \textsc{SubSegSievFac}, thus making $m|(n+j)/\Pi_j$ impossible.

We should also explain the purpose of the final loop in
\textsc{NewSegSievFac}. (The loop in the classical procedure
\textsc{SegSievFac} is identical, and plays the same role.)
Once we take care of all $m\leq \sqrt{n+\Delta}$ (and possibly some beyond),
what we have, for each $-\Delta \leq j\leq \Delta$, is either
that $\Pi_j = 1$, and $F_j$ contains a full factorization of $n+j$,
or $\Pi_j \ne 1$, and $F_j$ is missing a single large prime factor $p$,
which has to be equal to $n/\Pi_j$. We include that prime factor in the
factorization and are done.


\begin{algorithm}
  \caption{Main algorithm: 
factoring integers in $\lbrack n-\Delta,n+\Delta\rbrack$}\label{alg:segsievefac}
  \begin{algorithmic}[1]
    
    \Function{NewSegSievFac}{$n$,$\Delta$,$K$}
    \Ensure{for $-\Delta\leq j\leq \Delta$,
      $F_j$ is the list of pairs $(p,v_p(n+j))$ for $p|n+j$}
        \Require{$n,\Delta\in \mathbb{Z}^+$,   $\sqrt[3]{n}\leq \Delta<n$,
      $K\geq 5/2$}
        \State{$(F',\Pi') \gets \text{\textsc{SubSegSievFac}}(n-\Delta,2 \Delta,K \Delta)$}\Comment{Find factors $p\leq K\Delta$}
        \State{$F_j \gets F'_{j+\Delta}$, $\Pi_j \gets \Pi'_{j+\Delta}$
          for all $-\Delta \leq j\leq \Delta$}        
\State{$M\gets \lfloor K \Delta \rfloor + 1$}
\While{$M\leq \sqrt{n+\Delta}$}
\State{$R\gets \lfloor M \sqrt{\Delta/4 n}\rfloor$, $m_0 \gets M+R$}
\State{$\alpha_1 \gets \{- n/m_0^2\}$, $\alpha_0 \gets \{n/m_0\}$, 
$\eta\gets 5 \Delta/4 M$}
\State{$(a,a^{-1},q) \gets \text{\textsc{DiophAppr}}(\alpha_1,2 R)$}
\State{$c\gets \lfloor \alpha_0 q + 1/2\rfloor$,
  $k\gets \lfloor \eta q\rfloor$}
\For{$-k-1\leq j\leq k+1$}
\State{$r_0\gets - a^{-1} (c+j) \mo q$}
\For{$m\in (m_0 + r_0 + q \mathbb{Z}) \cap \lbrack M,M+2R\rbrack$}
\State{$n'\gets \lfloor (n+\Delta)/m\rfloor\cdot m$}  \Comment{$n'$ is a multiple of $m$}
\If{$n' \in \lbrack n-\Delta,n+\Delta\rbrack$}
\If{$m|(n'/\Pi_{n'-n_0})$} \Comment{$m$ is a new factor of $n'$}
\If{$m^2|n'$}
\State{$\Pi_{n'-n_0} \gets m^2\cdot \Pi_{n'-n_0}$, {\bf append} $(m,2)$ to $F_{n'-n}$}
\Else
\State{$\Pi_{n'-n_0} \gets m\cdot \Pi_{n'-n_0}$, {\bf append} $(m,1)$ to $F_{n'-n}$}
\EndIf
\EndIf
\EndIf
\EndFor
\EndFor
\State{$M\gets M+2 R+1$}
\EndWhile
\For{$n_0\leq n'\leq n_0+2 \Delta$}
\If{$\Pi_{n'-n_0} \ne n'$}
\State{$p_0\gets n'/\Pi_{n'-n_0}$, {\bf append} $(p_0,1)$ to $F_{n'-n_0}$}
\EndIf
\EndFor
\State{\Return{$F$}}
\EndFunction
  \end{algorithmic}
  \end{algorithm}

\section{Time analysis. Parameter choice.}

The space and time consumption of Algorithms
\ref{alg:erasieve}--\ref{alg:basterno}
is clearly as stated.


{\em Time consumption of main algorithm.}
Let us analyze the time consumption of Algorithm \ref{alg:segsieve}.
(Its space consumption, namely, $O(\Delta + \sqrt{K \Delta})$, will be clear.)
Sieving by primes $p\leq K\Delta$ gets done
by the traditional segmented-procedure
\textsc{SubSegSieve}$(n,2\Delta,K\Delta)$, which takes time
$O(K\Delta \log \log K \Delta)$ and space $O(\sqrt{K\Delta} +\Delta)$.
We must analyze now how much time it takes to sieve by integers
$K\Delta<m\leq \sqrt{n+\Delta}$.
(Our algorithm sieves by integers $K\Delta<m\leq \sqrt{n+\Delta}$,
not just by primes, simply because the Diophantine-approximation
algorithm cannot tell in advance which of the integers it outputs
will be prime. As we discussed before, we could apply a small wheel in the
hope of saving a constant factor in time.)

We will use the main ideas of Vinogradov's proof of the bound
\begin{equation}\label{eq:daremos}\sum_{n\leq x} \tau(n) = x \log x + O(x^{1/3} (\log x)^{5/3}),\end{equation}
as given in \cite[Ch.\ III, exer.\ 3--6]{MR0062138}\footnote{Actually,
  in the given reference, Vinogradov gives a bound of
  $O(x^{1/3} (\log x)^2)$ on the error term in (\ref{eq:daremos}).
  The reason is simply that he did not choose his parameters optimally:
  if the parameter $\tau$ in his proof is set to
  if the value of $\tau$ in the
$(A \log A)^{1/3}$ rather
  than $A$, the resulting bound is indeed as in (\ref{eq:daremos}).}, though
we will not use the bound (\ref{eq:daremos}) itself. Our treatment will
be self-contained.
  
Algorithm \ref{alg:segsieve} does not correspond
extremely closely to Vinogradov's approach -- 
we use our approximations $n/m_0$, $-n/m_0^2$ on the relatively broad intervals
on which they are useful, whereas Vinogradov's procedure changes 
approximations constantly. Nevertheless, we will be able to use the basic
approach that he took to bound an error term, though we of course will use
it to bound time consumption.

Let us look at an interval $\lbrack M,M+2R\rbrack$.
Finding $a$, $a^{-1}$, $q$ by Diophantine approximation (function
\textsc{DiophAppr}) takes, as we know, time $O(\log \max(Q,m_0^2)) =
O(\log n)$. The number of iterations of the
outer loop in \textsc{NewSegSieve} is
\begin{equation}\label{eq:seraph1}\ll \sqrt{\frac{4 n}{\Delta}}
\cdot \log \frac{\sqrt{n}}{K \Delta}
\leq \sqrt{\frac{n}{\Delta}} \log n
,\end{equation} and thus the total time taken
by the instructions inside the outer loop but outside the inner loop is
$O(\sqrt{n/\Delta} (\log n)^2)$.

The time it takes to go over and sieve by all $m\in \lbrack M,M+2R\rbrack$
congruent to $m_0 - a^{-1} (c+j) \mo q$ for all
$-k-1\leq j\leq k+1$ (where $k=\lfloor \eta q\rfloor$) is at most
\begin{equation}\label{eq:untu}
\begin{aligned}
  \left\lceil \frac{2 R+1}{q}\right\rceil \cdot (\lfloor \eta q\rfloor+ 3) &\leq
  \left(\frac{2 R}{q} + 2\right) (\eta q + 3) \ll
(R+q) \eta + \frac{R}{q}  + 1\\
&\ll \frac{\Delta^{3/2}}{\sqrt{n}} + \frac{R}{q} + 1,
\end{aligned}\end{equation}
since $R = \lfloor \sqrt{\Delta/4 n} \cdot M\rfloor$, 
$\eta = 5 \Delta/4 M$ and $q\leq R$.

Since the number of times the main loop is executed -- that is, the number
of intervals $\lbrack M,M+2R\rbrack$ we consider -- is
given by (\ref{eq:seraph1}),
the total contribution of the first and last terms in the last line
of (\ref{eq:untu}) is
\[\ll
\left(\frac{\Delta^{3/2}}{\sqrt{n}} + 1\right)
\sqrt{\frac{n}{\Delta}} \log n
\ll \Delta \log n,\]
since $\Delta\geq n^{1/3}$. It remains to account for the contribution
of $R/q$.
We may assume $M$ is large enough for $R$ to be $\geq 3$, as otherwise
the contribution of $R/q$ is bounded by thrice the contribution of the last
term $1$.

Now we proceed much as in \cite[Ch.\ III, exer.\ 3--6]{MR0062138}.
We will examine how $\alpha_1 = \{- n/m_0^2\}$ changes as $m_0$ increases;
we will then be able to tell how often Diophantine approximations $a/q$
to $\alpha_1$ with given $q$ can occur. To be precise:
we will see by how much $m_0$ has to increase for $-n/m_0^2$ to increase by
$1$ or more, and we also want to know for how long $-n/m_0^2$ can have
 a given, fixed Diophantine approximation $a/q$ as $m_0$ increases.

Consider two intervals $\lbrack m_0-R,m_0+R\rbrack$,
$\lbrack m_0'-R',m_0'+R'\rbrack$, where $m_0<m_0'$ and $R\leq R'$. Then
\begin{equation}\label{eq:skllet}
\left(-\frac{n}{(m_0')^2}\right) - \left(-\frac{n}{m_0^2}\right)
= n \frac{(m_0')^2-m_0^2}{m_0^2 (m_0')^2}.
\end{equation}

If $\alpha_1 = -n/m_0^2$ is $a/q + O^*(1/q R)$ and
$\alpha_1' = - n/(m_0')^2$ is $a/q + O^*(1/q R')$, it follows that
$n/m_0^2 - n/(m_0')^2$ is $O^*(1/q R + 1/q R') = O^*(2/q R)$.
Suppose that this is the case.

Since $\Delta<n$, we know that
$M' := m_0' - R' \leq \sqrt{n+\Delta}< \sqrt{2 n}$,
$R'\leq M' \sqrt{\Delta/4 n} < M'/2$
and $m_0' = M'+R' < (\sqrt{2}+1/2) \sqrt{n}< 2\sqrt{n}$; in the same way,
$m_0 = M+R < 3 M/2$. 
Clearly, $m_0'\leq 2 m_0$, as otherwise
$n/m_0^2 - n/(m_0')^2 \geq 3 n/4 m_0^2\geq 3 n/(m_0')^2 > 3/4$,
giving us a contradiction to $2/q R \leq 2/R \leq 2/3$.
Hence \[\frac{2}{q R}\geq
n \frac{(m_0')^2-m_0^2}{m_0^2 (m_0')^2} = n \cdot \frac{m_0'-m_0}{m_0^2}\cdot
\frac{m_0'+m_0}{(m_0')^2} \geq \frac{3 n}{4}\cdot \frac{m_0'-m_0}{m_0^3}.
\]
In other words, when the same approximant $a/q$ is valid at $m_0$ and $m_0'$,
\[m_0'-m_0 \leq \frac{8}{3} \frac{m_0^3}{q R n}.\]

The number of intervals for which $\alpha_1 = -n/m_0^2$ has a given
approximation $a/q$ is thus at most
\[\begin{aligned}\frac{8}{3} \frac{m_0^3}{q R (2 R + 1) n} + 1&\ll
 \frac{m_0^3}{q (R + 1)^2 n} + 1\leq
 \frac{m_0^3}{q M^2 \Delta/4} + 1\\
&< (3/2)^2 \frac{4 m_0}{q \Delta} + 1 \ll
\frac{m_0}{q \Delta} + 1.\end{aligned}\]

We should also see when $n/m_0^2-n/(m_0')^2 \geq 1$, or rather
when $n/m_0^2-n/(m_0')^2 \geq 1-2/R$, as then $n/m_0^2$ and $n/(m_0')^2$
may have Diophantine approximations that, while distinct, are congruent
$\mo 1$, i.e., differ by an integer. 
By (\ref{eq:skllet}), the inequality
$n/m_0^2-n/(m_0')^2 \geq 1-2/R$
is fulfilled exactly when
\[n ((m_0')^2-m_0^2) \geq m_0^2 (m_0')^2 (1-2/R),\]
and, since $R\geq 3$, that implies
\[6 n (m_0' - m_0) > m_0^3.\]
Hence, for given $m_0$, it makes sense to consider all following intervals with
$m_0'\leq m_0 + m_0^3/6 n$. In that range, $n/m_0^2$ and $n/(m_0')^2$ can
have the same Diophantine approximation $\mo 1$ only if they in fact have the
same Diophantine approximation.

Since $R$ increases and the intervals are of width $2 R$,
there will be at most
$m_0^3/12 R n + 1 \ll m_0^2/\sqrt{\Delta n} + 1\ll m_0^2/\sqrt{\Delta n}$
intervals with $m_0'\leq m_0+m_0^3/6 n$.
(Note that $m_0^2/\sqrt{\Delta n} \geq K^2 \Delta^2/\sqrt{\Delta n}> 1$.)
Among those intervals,
$\ll m_0'/q \Delta + 1 \ll m_0/q \Delta + 1$ will have a given approximation $a/q \mo 1$.

As we said before, we have to account for the total contribution of $R/q$.
Since $1/q$ and $1/q^2$ decrease as $q$ increases, the worst-case scenario
is for there to be as many $m_0$'s as possible for which the approximation
has $q$ as small as possible. We would have
$\lfloor O(m_0/q \Delta+1) \phi(q) \rfloor = O(m_0/\Delta+q)$ values of $a/q$ with $q\leq Q$, and none for
$q>Q$, where $Q\ll m_0/(\Delta n)^{1/4}$.
(To bound $Q$, we apply
$\sum_{q\leq Q} \phi(q) = |\{1\leq q_1,q_2\leq Q: \text{$q_1$, $q_2$ coprime}\}| \gg Q^2$.)
The contribution of $R/q$ for
all intervals
with $m_0'\leq m_0+m_0^3/6 n$ is thus 
\[\ll R \sum_{q\leq Q} \left(\frac{m_0}{\Delta q} + 1\right)
\ll R \cdot \left(\frac{m_0}{\Delta} \log Q + Q\right).\]
Now split $\lbrack M_0, 2 M_0\rbrack$ into
$O(n/M_0^2)$ chunks of the form $\lbrack m_0, m_0 (1+ m_0^2/6 n)\rbrack$.
We obtain that the contribution for all intervals with $m_0$ inside
$\lbrack M_0,2 M_0\rbrack$ is
\[\begin{aligned} &\ll \frac{n}{M_0^2} \cdot 
2 M_0 \sqrt{\frac{\Delta}{4 n}} \cdot  \left(\frac{M_0}{\Delta} \log n +
\frac{M_0}{(\Delta n)^{1/4}}\right)\\
&\ll \frac{\sqrt{n}}{\sqrt{\Delta}} \log n +(n \Delta)^{1/4},\end{aligned}\]
and thus the total contribution will be
\[\ll
\frac{\sqrt{n}}{\sqrt{\Delta}} (\log n)^2 +(n \Delta)^{1/4} \log n.\]

We conclude that the total time consumption of Algorithm
\ref{alg:segsieve} is
\[\begin{aligned}
&\ll \Delta \log \log K \Delta + \Delta \log n + 
\frac{\sqrt{n}}{\sqrt{\Delta}} (\log n)^2 +(n \Delta)^{1/4} \log n\\
&\ll \Delta \log n + \sqrt{\frac{n}{\Delta}} (\log n)^2,\end{aligned}\]
since $\Delta\geq n^{1/3}$. Under the stronger assumption
$\Delta \geq n^{1/3} \log^{2/3} n$, we obtain that the total time consumption
is
\[O(\Delta \log n).\]


{\em Time consumption
  of Algorithm \ref{alg:segsievefac}.} We must now analyze
Algorithm \ref{alg:segsievefac},
which factors integers
in the interval $\lbrack n-\Delta,n+\Delta\rbrack$. 
Everything is much the same as for Algorithm \ref{alg:segsieve},
except in one respect. The total space taken
will be
\[O(\Delta \log n)\]
rather than $O(\Delta)$, simply because storing the list of prime factors
of an integer in $\lbrack n-\Delta,n+\Delta\rbrack$
takes $O(\log n)$ bits. (To wit: the number of bits taken to store a factor
$p^\alpha$ is $O(\log p) + O(\log(\alpha+1))$, which is $O(\alpha \log p)$.
Hence, storing all the factors $p_i^{\alpha_i}$ of an integer
$n = p_1^{\alpha_1}\dotsb p_k^{\alpha_k}$ takes
$O\left(\sum_i \alpha_i \log p_i\right) = O(\log n)$ bits.)


The total time estimate is still
\[O\left(\Delta \log n\right).\]
The claims in (\ref{eq:rada1}) and
(\ref{eq:rada2}) follow immediately once one splits the interval
$\lbrack 1,N\rbrack$ into intervals of length $2\Delta$.
The main theorem is thus proved.

\section{Further perspectives}
   
   It is tempting to try to improve on this algorithm by taking a longer truncated Taylor
   expansion in (\ref{eq:anzio}). However, this would require us to
   find the solutions $x\in\{-R,-R+1,\dotsc,R\}$ to 
   a $P(x)\in \lbrack -\eta,\eta\rbrack \mo \mathbb{Z}$,
   where $P$ is a polynomial with
   $\deg P \geq 2$ and $\eta\ll 1/R$ or thereabouts.

Solving a quadratic modular equation
$a_2 x^2 + a_1 x + a_0 \equiv 0 \mo q$ is not the main difficulty in our
context. We can obviously reduce such a problem to taking square-roots
$\mo q$. Now, the problem of finding square-roots modulo $q$
(for $q$ arbitrary) is well-known to be equivalent to
   factoring $q$ \cite{Rabin:1979:DSP:889813}.
(There is a classical algorithm (Tonelli-Shanks) for finding square-roots 
to prime modulus.) This is not a problem, since we can factor
all integers $q\leq x^{1/4}$ (say) in advance, before we start sieving.

The problem is that it is not at all clear how to reduce 
finding solutions to $P(x)\in \lbrack -\eta,\eta \rbrack \mo \mathbb{Z}$,
$\deg P=2$,
to finding solutions to quadratic equations $\mo q$. We can try to find
rational approximations with the same denominator $q$ ({\em simultaneous
Diophantine approximation}) to the non-constant coefficients of $P$,
but such approximations will be generally worse than when we approximate a single
real number, and so matters do not work out: we need a more, not less,
precise approximation than before to the leading coefficient, since it
is now the coefficient of a quadratic term.

Concretely, for $P(x) = \alpha_2 x^2 + \alpha_1 x + \alpha_0$,
in order to reduce the problem of finding solutions to
$P(r) \in \lbrack -\eta,\eta \rbrack \mo \mathbb{Z}$ with
$r\in \lbrack R,R\rbrack \cap \mathbb{Z}$ to the problem of solving
a quadratic modular equation,
 we would need
$a_1$, $a_2$, $q$ with $|\alpha_2 - a_2/q|\leq 1/q R^2$,
$|\alpha_1 - a_1/q|\leq 1/q R$. In general, we can do such a thing only for 
$q\gg R^3$, and that is much too large. For one thing, for $\epsilon\sim 1/R$,
we would have to solve $\epsilon q \sim R^2$ distinct equations, and that is
obviously too many.

   \subsection{Geometric interpretation}
As we have seen, sieving integers in the interval 
$\lbrack n-\Delta,n+\Delta\rbrack$ reduces to finding integers $m$ such
that $\{n/m\}$ is close to $0$ modulo $1$.
This task is equivalent to finding integer points
close to a hyperbola $x\mapsto n/x$. Seen from this perspective, our
approach consists simply in approximating a hyperbola locally by linear
functions. Such a geometric perspective is already present (and dominant) in
 \cite{zbMATH02656343}.

What we did in \S \ref{sec:analys} then amounts to finding points close
to a line, starting with an approximation of the slope by a rational, and then
proceeding by modular arithmetic. 

Taking one more term in the Taylor expansion would be
 the same as approximating a hyperbola by segments of parabolas, rather than 
by segments of lines. We could then take longer segments, and thus
hope to capture points near the hyperbola efficiently even when 
$\Delta$ is considerably smaller than $n^{1/3}$. (We can hope to
capture them efficiently because the segments are long enough that
we expect at least one point in each segment.) The problem of how to find
the points remains. It seems difficult to do so without reducing a
non-linear problem $\mo 1$ to a problem $\mo q$, and we do not yet know how to carry out such a reduction
well.

\section{A few words on the implementation}\label{app:volovoi}

I have written and tested a simple 
implementation as proof-of-purpose, that is, mainly so as to check that the algorithms work as described.
On the higher end of what can be done in 64-bit
computer arithmetic (say: $n=5\cdot 10^{18}$ and $\Delta = 2\cdot 10^7$ or
$\Delta = 4\cdot 10^7$),
my implementation of algorithm \textsc{NewSegSiev}
runs substantially faster than my own implementation of the traditional
algorithm \textsc{SegSiev}. Still, on those same inputs, my implementation of
\textsc{NewSegSiev} is clearly slower (by a factor between $2$ and
$2.5$) than a publicly
available, highly optimized implementation \cite{primesieve} of the traditional
algorithm (basically \textsc{SegSiev}, but improved in all the ways
detailed below, and some other ones) on the same interval.
Diophantine approximation (Algorithm \ref{alg:diophfrac})
turns out to take less than $2$ percent of
the running time of \textsc{NewSegSiev}, so
the fact that \cite{primesieve} runs faster is due to better
coding, and not to the overhead of Diophantine approximation. For the values of
$n$ and $\Delta$ above, we are talking about
total running times of only a few seconds in all cases.

The advantage of \textsc{NewSegSieve} over existing methods
should become clearer as $n$ grows larger --
meaning substantially larger than $2^{64}$. 
However, on the range $n>2^{64}$,
it seems harder to find highly optimized implementations
of the traditional algorithm for comparison.

Here are a few hints for the reader who would like to write a more serious
program on his or her own. Most of these tricks are standard, but are scattered
here and there in the literature (and in code).

\begin{enumerate}
\item Obviously, 
  we can save on space by storing the sieve as a bit array. Saving on space
  leads to better cache usage, and hence often to time savings as well.
  Our very conventions for space usage (expressed in terms of bits, rather
  than words) reflect this fact. Of course, in some ranges, it can save time
  to be a little more wasteful and store some integer data as words
  (e.g., prime factors, in the case of the sieve applied to factorization).
\item We can first apply a simple sieve to the integers between
  $1$ and $M = \prod_{p\leq p_0} p$ (where $p_0=17$, say), taking only the
  effect of primes $p\leq p_0$ into account (whether we are sieving for
  primality, or computing the M\"obius function
  $\mu$, or factoring numbers:  in the $n$th entry, we would store whether $n$
  is coprime to $P$, or, instead, store $\mu(\gcd(n,M^2))$, if we
  are computing $\mu$,
  or the set 
  $\{p\leq p_0: p|n\}$, if we are factoring numbers). We then initialize
  our sieve by repeating that block of length $M$, and so we do not need
  to sieve by the primes $p\leq p_0$ ever again.

\item We can of course implement a sieve in parallel in a trivial sense,
  by letting different processors look at different intervals (of length
  at least $n^3 (\log n)^{2/3}$, in our case). There seems to be a small
  literature on parallel implementations of sieves; see, e.g.,
  \cite{MR1294306}.

\item At the end of \S \ref{sec:analys}, we went briefly over the issue of
  ``false alarms'', that is, values of $m$ that lie in a valid congruence
  class $m\equiv m_0+r_0 \mo q$ but do not actually have
  multiples within the interval $I$ we are sieving. Such false alarms do not
  change the outcome of the algorithm, but they do waste
  some time.  It is possible to avoid them; the details are given in the
  first version of the present paper, available on arXiv.org. There, in fact,
  the algorithm without ``false alarms'' is called \textsc{NewSegSiev},
  and presented as the main version of
  the algorithm.
  Unfortunately, at least in my implementation, the process of eliminating
  false alarms takes more time than it saves.
  
\item As we said in the introduction, Oliveira e Silva has shown
  how to use cache efficiently when implementing a sieve of Eratosthenes
  \cite{OeS}, \cite[Algorithm 1.2]{OSHP}. In effect, when sieving an interval
  of length $L$, he needs not much more than $\sqrt{L}$ units of memory in
  cache. There seems to be no reason why Oliveira e Silva's technique
  can't be combined with the algorithms put forward here. Thus we could
  hope to sieve intervals of the form $\lbrack n-\Delta, n+\Delta\rbrack$,
  $\Delta\sim n^{1/3} (\log n)^{2/3}$, while using no more than
  $O(n^{1/6} (\log n)^{1/3})$ units of memory in cache at a time (and
  $O(n^{1/3} (\log n)^{2/3})$ units of memory in total). The range up to about
  $n\sim 10^{36}$ could then become accessible, at least for short or medium-sized intervals. Note that we can still use
  128-bit-integer arithmetic in that range.
  \end{enumerate}
\bibliographystyle{alpha}
\bibliography{sieve}
\end{document}